\documentclass[12pt,oneside,english]{amsart}
\usepackage[latin1]{inputenc}
\usepackage{amssymb}

\makeatletter
\newtheorem{theorem}{Theorem}

\usepackage{babel}

\makeatother
\begin{document}
\baselineskip=17pt

\title[On connection between the numbers of permutations]{On connection between the numbers of permutations and full cycles with some restrictions on positions and up-down structure}

\author{Vladimir Shevelev}
\address{Departments of Mathematics \\Ben-Gurion University of the
 Negev\\Beer-Sheva 84105, Israel. e-mail:shevelev@bgu.ac.il}


\begin{abstract}
We discuss both simple and more subtle connections between the numbers of permutations and full cycles with some restrictions, in particular, between the numbers of permutations and full cycles with prescribed up-down structure.
\end{abstract}

\maketitle

\section{Introduction}
The following theorem is very well known.

\begin{theorem}\label{t1}  \upshape (cf, e.g.,[6]).\slshape \enskip The number of full cycles of $n$ elements
$1,2,\ldots, n $ equals to $(n-1)!$ i.e. the number of all permutations of $n-1$ elements
$1,2,\ldots, n-1 $.
\end{theorem}
     Let A be a quadratic (0,1)-matrix of order $n$. For a permutation $\sigma$ of elements $1,2,\ldots, n $ denote $H_\sigma$ the incidence (0,1)-matrix of $\sigma$. Let us consider the set $B(A)$ of permutations $\sigma$ for which $H_\sigma\leq A$.
$B(A)$ is the class of permutations with restriction on positions which is defined by zeros of the matrix $A$.   It is well known that $|B(A)|=per A$. Furthermore, let us  consider a matrix function which is defined, similar to permanent, by the formula

\begin{equation}\label{1}
perf  A= \sum_{f.c.\sigma}\sum^n_{i=1}a_{i,\sigma(i)},
\end{equation}

where the external sum is over all full cycles of elements $1,2,\ldots, n$. In particular,

$$
perf (a_{11})=a_{11},
$$

$$
perf\begin{pmatrix}a_{11}\;\; a_{12}\\ a_{21}\;\; a_{22} \end{pmatrix}=a_{12}\;a_{21},
$$
\newpage

$$
perf \begin{pmatrix}a_{11}\;\; a_{12}\;\;a_{13}\\ a_{21}\;\; a_{22}\;\;a_{23}\\a_{31}\;\;a_{32}\;\;a_{33} \end{pmatrix}=a_{12}\;a_{23}\;a_{31}+a_{31}\;a_{21}\;a_{32},
$$

$$
perf\begin{pmatrix}a_{11}\;\; a_{12}\;\;a_{13}\;\;a_{14}\\ a_{21}\;\; a_{22}\;\;a_{23}\;\;a_{24}\\a_{31}\;\;a_{32}\;\;a_{33}\;\;a_{34}\\ a_{41}\;\; a_{42}\;\;a_{43}\;\;a_{44} \end{pmatrix}=a_{11}\;a_{23}\;a_{34}\;a_{41}+a_{12}\;a_{24}\;a_{31}\;a_{43}+
$$

$$
+a_{13}\;a_{23}\;a_{34}\;a_{42}+a_{13}\;a_{24}\;a_{32}\;a_{41}+
a_{14}\;a_{23}\;a_{31}\;a_{42}+a_{14}\;a_{21}\;a_{32}\;a_{43}.
$$

From our results \cite{7}, where we considered a system of "partial permanents" including
$perf A $ with notation $pper _1A,$  it follows an expansion of  $perf A$ by the first row of $A$.

\begin{theorem}\label{t2}  $perf A=\sum^n_{j=2}a_{1j}perf A^*_{1,j}$,  where $A^*_{1,j}$ is obtained from the minor matrix $A_{1,j}$  by a cyclic permutation of its first $j-1$ columns by the rule:  $$2\rightarrow 1,\;\;3\rightarrow 2 ,\ldots, j-1\rightarrow j-2,\enskip 1\rightarrow j-1.$$

\end{theorem}

For example,

$$
perf \begin{pmatrix}a_{11}\;\; a_{12}\;\;a_{13}\;\;a_{14}\\ a_{21}\;\; a_{22}\;\;a_{23}\;\;a_{24}\\a_{31}\;\;a_{32}\;\;a_{33}\;\;a_{34}\\ a_{41}\;\; a_{42}\;\;a_{43}\;\;a_{44} \end{pmatrix}=a_{12}perf \begin{pmatrix}a_{21}\;\;a_{23}\;\;a_{24}\\
a_{31}\;\;a_{33}\;\;a_{34}\\a_{41}\;\;a_{43}\;\;a_{44}\end{pmatrix}+
$$

$$
+a_{13}perf \begin{pmatrix}a_{22}\;\; a_{21}\;\;a_{24}\\ a_{32}\;\; a_{31}\;\;a_{34}\\a_{42}\;\;a_{41}\;\;a_{44} \end{pmatrix}+a_{14}perf\begin{pmatrix}a_{22}\;\;a_{23}\;\;a_{21}\\
a_{32}\;\;a_{33}\;\;a_{31}\\a_{42}\;\;a_{43}\;\;a_{41}\end{pmatrix}.
$$

The algorithm of Theorem 2 is easily realized by computer. The observations obtained using this algorithm allowed us to formulate the following conjecture (1993) which is a generalization of Theorem 1.

\bfseries Conjecture  1\mdseries \cite{8}. \slshape Let $P^{(i)}_n$ be quadratic $(0,1)$ matrix of order $n$ with only 1`s on places $(1,1+i), (2,2+i),\ldots,(n-i,n),\;0\leq i \leq n-1,$ $J_n $ be $n\times n$ matrix composed of 1`s. Then \upshape

\begin{equation}\label{2}
perf(J_n-\sum^k_{j=1}P_n^{(j)})=per(J_{n-1}-\sum^{k-1}_{i=0}P^{(i)}_{n-1}).
\end{equation}

It is clear that, in the case of $k=0,$ we obtain Theorem 1 in the form
\newpage

\begin{equation}\label{3}
perf(J_n)=per(J_{n-1})=(n-1)!
\end{equation}

In 1994, this conjecture was proved independently by Ira M.Gessel using ideals of his paper \cite{1} and Richard P.Stanley which gave a direct proof (private correspondences, unpublished).

     In this paper we discuss  quite another intrigueing connections between the numbers
of permutations and full cycles with prescribed up-down structure.

\section{On up-down basis polynomials}

     \slshape Basis polynomial with up-down index $k$\upshape, denoted by $\left\{\begin{matrix} n\\
k\end{matrix}\right\},$ is the number of permutations $\pi=(\pi_1, \pi_2,\hdots, \pi_n)$ of elements $1,2,\hdots, n$
with the following condition: for $i\in [1,n-1],$ we have $\pi_i< \pi_{i+1}$ if and only if in the (n-1)-digit binary expansion of $k$ the i-th digit equals to 1 \cite{10}.

Let $k\in[2^{t-1}, 2^t)$ and the $(n-1)$-digit binary expansion of
$k$ have a form:

\begin{equation}\label{4}
k=\underbrace{0 \ldots 0}_{n-t-1}\; 1 \;\underbrace{0 \ldots
0}_{s_2-s_1-1}\; 1 \;\underbrace{0 \ldots 0}_{s_3-s_2-1}\; 1\ldots 1
\;\underbrace{0 \ldots 0}_{s_m-s_{m-1}-1}\; 1 \;\underbrace{0 \ldots
0}_{t-s_m}\enskip,
\end{equation}

where

$$
1=s_1< s_2 <\ldots < s_m
$$
are places of 1's \slshape after \upshape $n-t-1$ 0's before the
first 1.

In \cite{9}, using the fundamental Niven's result \cite{5}, the following formula was proved.

\begin{theorem}\label{t3}
$$
\left\{\begin{matrix} n \\ k
\end{matrix}\right\}=(-1)^m+$$
\begin{equation}\label{5}
\sum^m_{p=1}(-1)^{m-p}\sum_{1\leq i,< i_2< \hdots< i_p\leq m}
\left( \begin{matrix} n \\ t+1-s_{i_p}\end{matrix}\right)\prod
^p_{r=2}\begin{pmatrix} n-t+s_{i_r}-1
\\ s_{i_r}-s_{i_{r-1}}\end{pmatrix}.
\end{equation}
\end{theorem}

 Let us write (\ref{4}) in the form

\begin{equation}\label{6}
k=2^{t_1-1}+2^{t_2-1} +\hdots + 2^{t_m-1} ,\;\;t_1> t_2 > \hdots > t_m\geq 1.
\end{equation}

Comparing (\ref{4}) and (\ref{6}), we find

\begin{equation}\label{7}
t_i=t-s_i+1,\;\;i=1,2,\hdots, m.
\end{equation}
\newpage
It is easy to check directly the following identity

\begin{equation}\label{8}
\begin{pmatrix} n\\t_{i_p}\end{pmatrix} \prod^p_{r=2}\begin{pmatrix} n-t_{i_r}\\ t_{i_{r-1}}-t_{i_r}\end{pmatrix}=
\begin{pmatrix} n\\ t_{i_1}\end{pmatrix}\begin{pmatrix} t_{i_1}\\ t_{i_2}\end{pmatrix}
\begin{pmatrix} t_{i_2}\\ t_{i_3}\end{pmatrix}\hdots\begin{pmatrix} t_{i_{p-1}}\\ t_{i_p}\end{pmatrix}.
\end{equation}

Now by (\ref{5}), (\ref{7}) and (\ref{8}) we obtain $\left\{\begin{matrix} n\\
k\end{matrix}\right\}$ as a linear combinations of binomial coefficients.

\begin{theorem}\label{t4} $($\cite{10}$)$.  For $k$ $(\ref{6})$ we have
\begin{equation}\label{9}
\left\{\begin{matrix} n\\
k\end{matrix}\right\}=(-1)^m+\sum^m_{p=1}(-1)^{m-p}\sum_{1\leq i_1 < i_2 < \ldots <i_p \leq m}\begin{pmatrix} n\\t_{i_1}\end{pmatrix}\prod^p_{r=2}\begin{pmatrix}t_{i_{r-1}}\\t_{i_r}\end{pmatrix}.
\end{equation}
\end{theorem}

     Below, as in \cite{10} we consider $\left\{\begin{matrix} n\\
k\end{matrix}\right\}$ from the formal (wider than only combinatorial) point of view: according to (\ref{9}) it is a polynomial in $n$ of degree $t_1=\lfloor\log_2(2k)\rfloor$. In particular,

\begin{equation}\label{10}
\left\{\begin{matrix} 0\\
k\end{matrix}\right\}=(-1)^m=\tau_k,
\end{equation}

where $ \tau_k,\;\; k=0,1,2, \ldots $, is the Thue-Morse sequence \cite{4}, \cite{2}.

The following theorem is equivalent to Theorem 10 \cite{10}. Here we give a more detailed proof of this theorem.

\begin{theorem}\label{t5}
For $k$ $(\ref{6})$ we have

\begin{equation}\label{11}
\left\{\begin{matrix} n\\
k\end{matrix}\right\}=\left |\begin{matrix}\begin{pmatrix} n\\ t_1 \end{pmatrix}&\begin{pmatrix} n\\ t_2 \end{pmatrix}&
\begin{pmatrix} n\\ t_3 \end{pmatrix}& \ldots &\begin{pmatrix} n\\ t_{m-1} \end{pmatrix}&\begin{pmatrix} n\\ t_m \end{pmatrix}& 1\\ 1 & \begin{pmatrix} t_1\\t_2\end{pmatrix} & \begin{pmatrix} t_1\\t_3\end{pmatrix}&\ldots & \begin{pmatrix} t_1\\t_{m-1}\end{pmatrix}& \begin{pmatrix} t_1\\t_m \end{pmatrix}& 1 \\ 0 & 1 & \begin{pmatrix} t_2\\t_3\end{pmatrix}& \ldots & \begin{pmatrix} t_2\\t_{m-1}\end{pmatrix}& \begin{pmatrix} t_2\\t_m\end{pmatrix} & 1\\
0 & 0 & 1 & \ldots & \begin{pmatrix} t_3\\ t_{m-1}\end{pmatrix} &\begin{pmatrix} t_3\\ t_m \end{pmatrix} & 1\\ \ldots & \ldots & \ldots & \ldots & \ldots & \ldots\\ 0 & 0 & 0 & \ldots & 1 & \begin{pmatrix} t_{m-1}\\ t_m \end{pmatrix} & 1\\ 0 & 0 & 0 & \ldots & 0 & 1 & 1 \end{matrix}\right|.
\end{equation}
\end{theorem}
\newpage

\bfseries Proof.\mdseries  The number of diagonals of matrix (\ref{11}) having no 0's equals to permanent of the following $(m+1)\times (m+1)$ matrix

$$
C_{m+1}=\begin{pmatrix} 1 &1 &1 &\ldots &1 &1 &1\\1 &1 &1 &\ldots &1 &1 &1\\0 &1 &1 &\ldots &1 &1 &1\\
0 &0 &1 &\ldots &1 &1 &1\\\ldots &\ldots &\ldots &\ldots &\ldots &\ldots &\ldots \\0 &0 &0 &\ldots &1 &1 &1\\
0 &0 &0 &\ldots &0 &1 &1\end{pmatrix}.
$$

Decomposing $ per C_{m+1}$ by the last row we find

$$
per C_{m+1}=2 per C_m = 2^2 per C_{m-1}=\ldots =2^m.
$$

Denote $A\;(m+1)\times(m+1)$ matrix in (\ref{11}) and consider $(m\times m)$ upper-triangle submatrix $T$ with the main diagonal composed of 1's. Let us choose $p$ 1's of the main diagonal of $T$ in its rows $(1\leq)i_1 < i_2 < \ldots < i_p \leq m$. To such choice corresponds a diagonal of $A$ composed of the other $(m-p)$ 1's of the main diagonal of T + the unit in the last column of A which is the continuation of the $i_p$-th row of $T+ elements \begin{pmatrix} n\\ t_{i_1}\end{pmatrix},\begin{pmatrix} t_{i_1}\\ t_{i_2}\end{pmatrix}, \begin{pmatrix} t_{i_2}\\ t_{i_3}\end{pmatrix},\ldots, \begin{pmatrix} t_{i_{p-1}}\\ t_{i_p}\end{pmatrix}$ such that in all we have $m-p+1+p=m+1$ elements of A which are in different rows and columns. As a result, we obtain $\sum^m _{p=0}\begin{pmatrix} m\\ p \end {pmatrix}=2^m,$ i.e. \slshape all \upshape diagonals of A having no 0's (note that, to $p=0$ corresponds the choice of the emply subset of 1's of the main diagonal of T, i.e. all these 1's and the unit in the first row of A form in this case the only diagonal of 1's). Therefore,

$$
per A= 1+\sum^m_{p=1}\sum_{1\leq i_1 < i_2 < \ldots < i_p\leq m}\begin{pmatrix} n\\ t_{i_1}\end{pmatrix}
\prod^p_{r=2}\begin{pmatrix} t_{i_{t-1}}\\t_{i_r}\end{pmatrix}.
$$
What is left-to find the number of cycles of diagonal which contains elements
$$\begin{pmatrix} n\\ t_{i_1}\end{pmatrix},\begin{pmatrix} t_{i_1}\\ t_{i_2}\end{pmatrix}, \begin{pmatrix} t_{i_2}\\ t_{i_3}\end{pmatrix},\ldots, \begin{pmatrix} t_{i_{p-1}}\\ t_{i_p}\end{pmatrix}$$
and 1 in position $(i_p+1,\enskip m+1).$ Put, for the symmetry, $n=t_0.$ Then matrix in (11) becomes "quasi-Toeplitz"
in the sense that, for a fixed $j\geq0,$ the diagonal $(a_{i,\enskip i+j})_{i\geq1}$ contains only non-zero elements of the form:
$$ \begin{pmatrix} t_0\\ t_{i_{j+1}}\end{pmatrix}, \enskip \begin{pmatrix} t_1\\ t_{i_{j+2}}\end{pmatrix},...$$
with the fixed difference of subindices, which is equal to $j+1.$ Note that,\newpage  every non-zero element of the form $(a_{i,\enskip i+j})$ is contained in a cycle of the form $[(a_{i,\enskip i+j}), 1,...,1]$ with $j$ 1's in the positions $ (2,1), \enskip (3,2),..., (j+1,j),$ i.e. it is contained in a cycle of length $j+1.$ Thus the considered elements are in cycles of length $i_1,\enskip i_2-i_1,...,i_p-i_{p-1},\enskip m+1-i_p$ ( the last length corresponds to 1 in position $(i_p+1,\enskip m+1)).$ Show that we have $p+1$ distinct cycles. Indeed, two cycles either coincide or disjoint. The total length of all cycles is
$$ i_1+(i_2-i_1)+...+(i_p-i_{p-1})+m+1-i_p=m+1.$$
If some two cycles coincide, then these two cycles are considered as one and, in this case, the total sum will be less than $m+1,$ which is impossible for a diagonal. Thus we have exactly $p+1$ cycles of the considered diagonal and the decrement of the corresponding permutation equals to $m+1-(p+1)=m-p$ and the parity of this permutation is $(-1)^{m-p}.$ This completes proof of formula (\ref{9}).$\blacksquare$
 \newline
\indent Many different properties of $\left\{\begin{matrix} n\\k \end{matrix}\right\}$ were proved in [10]. Let us prove an additional interesting property.

\begin{theorem}\label{t6}  $($cf.our comment to sequence A060351 \cite{11}$)$.

If $k\equiv 0\pmod 4,$ then
\begin{equation}\label{12}
\left\{\begin{matrix} n\\k \end{matrix}\right\}-\left\{\begin{matrix} n\\k+1 \end{matrix}\right\}+
\left\{\begin{matrix} n\\k+2\end{matrix}\right \}-\left\{\begin{matrix} n\\k+3 \end{matrix}\right\}=0.
\end{equation}
\end{theorem}

\bfseries Proof. \mdseries  By the condition,

$$
k=2^{t_1-1}+2^{t_2-1}+ \ldots +2^{t_m-1}.\;\; t_1>t_2>\ldots >t_m\geq 3
$$

and, by (\ref{11}), we have

$$
\left\{\begin{matrix} n\\k +3\end{matrix}\right\}=\left|\begin{matrix}\begin{pmatrix} n\\ t_1 \end{pmatrix} &
\begin{pmatrix} n\\ t_2 \end{pmatrix} & \ldots & \begin{pmatrix} n\\ t_m \end{pmatrix} &
\begin{pmatrix} n\\ 2 \end{pmatrix} & \begin{pmatrix} n\\ 1 \end{pmatrix} & 1\\
1 & \begin{pmatrix} t_1\\ t_2 \end{pmatrix}& \ldots & \begin{pmatrix} t_1\\ t_m \end{pmatrix}&
\begin{pmatrix} t_1\\ 2\end{pmatrix}& \begin{pmatrix} t_1\\1 \end{pmatrix} & 1\\
0 & 1 & \ldots & \begin{pmatrix} t_2\\ t_m \end{pmatrix} & \begin{pmatrix} t_2\\ 2 \end{pmatrix}&
\begin{pmatrix} t_2\\ 1 \end{pmatrix} & 1  \\
\ldots &\ldots & \ldots & \ldots & \ldots & \ldots & \ldots\\
0 & 0 & \ldots & \begin{pmatrix} t_{m-1}\\ t_m \end{pmatrix} &  \begin{pmatrix} t_{m-1}\\ 2 \end{pmatrix}&
\begin{pmatrix} t_{m-1}\\ 1\end{pmatrix} & 1\\
0 & 0 & \ldots & 1 &  \begin{pmatrix} t_m\\ 2 \end{pmatrix}&
\begin{pmatrix} t_m\\ 1\end{pmatrix} &  1\\
0 & 0 & \ldots & 0 & 1 & 2 & 1\\ 0 & 0 & \ldots & 0 & 0 & 1 & 1
\end{matrix}\right|=
$$\newpage

$$
 - \left\{\begin{matrix} n\\k +2\end{matrix}\right\}+\left|\begin{matrix}\begin{pmatrix} n\\ t_1 \end{pmatrix} &
\begin{pmatrix} n\\ t_2 \end{pmatrix} & \ldots & \begin{pmatrix} n\\ t_m \end{pmatrix} &
\begin{pmatrix} n\\ 2 \end{pmatrix} & \begin{pmatrix} n\\ 1 \end{pmatrix} \\
1 & \begin{pmatrix} t_1\\ t_2 \end{pmatrix}& \ldots & \begin{pmatrix} t_1\\ t_m \end{pmatrix}&
\begin{pmatrix} t_1\\ 2\end{pmatrix}& \begin{pmatrix} t_1\\1 \end{pmatrix} \\
0 & 1 & \ldots & \begin{pmatrix} t_2\\ t_m \end{pmatrix} & \begin{pmatrix} t_2\\ 2 \end{pmatrix}&
\begin{pmatrix} t_2\\ 1 \end{pmatrix}  \\
\ldots &\ldots & \ldots & \ldots & \ldots & \ldots \\
0 & 0 & \ldots & \begin{pmatrix} t_{m-1}\\ t_m \end{pmatrix} &  \begin{pmatrix} t_{m-1}\\ 2 \end{pmatrix}&
\begin{pmatrix} t_{m-1}\\ 1\end{pmatrix} \\
0 & 0 & \ldots & 1 &  \begin{pmatrix} t_m\\ 2 \end{pmatrix}&
\begin{pmatrix} t_m\\ 1\end{pmatrix} \\
0 & 0 & \ldots & 0 & 1 & 2
\end{matrix}\right|=
$$

$$
- \left\{\begin{matrix} n\\k +2\end{matrix}\right\}-\left|\begin{matrix}\begin{pmatrix} n\\ t_1 \end{pmatrix} &
\begin{pmatrix} n\\ t_2 \end{pmatrix} & \ldots & \begin{pmatrix} n\\ t_{m-1} \end{pmatrix} &
\begin{pmatrix} n\\ t_m \end{pmatrix} & \begin{pmatrix} n\\ 1 \end{pmatrix} \\
1 & \begin{pmatrix} t_1\\ t_2 \end{pmatrix}& \ldots & \begin{pmatrix} t_1\\ t_{m-1} \end{pmatrix}&
\begin{pmatrix} t_1\\t_m\end{pmatrix}& \begin{pmatrix} t_1\\1 \end{pmatrix} \\
0 & 1 & \ldots & \begin{pmatrix} t_2\\ t_{m-1} \end{pmatrix} & \begin{pmatrix} t_2\\t_m \end{pmatrix}&
\begin{pmatrix} t_2\\ 1 \end{pmatrix}  \\
\ldots &\ldots & \ldots & \ldots & \ldots & \ldots \\
0 & 0 & \ldots & 1 &\begin{pmatrix} t_{m-1}\\ t_m \end{pmatrix} &  \begin{pmatrix} t_{m-1}\\ 1 \end{pmatrix}
\\
0 & 0 & \ldots & 0 &  1&
\begin{pmatrix} t_m\\ 1\end{pmatrix}
\end{matrix}\right|+
$$

\begin{equation}\label{13}
2\left|\begin{matrix} \begin{pmatrix} n\\ t_1 \end{pmatrix}&\begin{pmatrix} n\\ t_2 \end{pmatrix}&
\ldots & \begin{pmatrix} n\\ t_{m-1} \end{pmatrix}&\begin{pmatrix} n\\ t_m \end{pmatrix}&\begin{pmatrix} n\\ 2 \end{pmatrix} \\
1 & \begin{pmatrix} t_1\\ t_2\end{pmatrix}& \ldots &\begin{pmatrix} t_1\\ t_{m-1 }\end{pmatrix}& \begin{pmatrix} t_1\\ t_m \end{pmatrix}& \begin{pmatrix} t_1\\ 2 \end{pmatrix}\\ 0 & 1 & \ldots &
\begin{pmatrix} t_2\\ t_{m-1} \end{pmatrix}& \begin{pmatrix} t_2\\ t_m \end{pmatrix}&\begin{pmatrix} t_2\\2 \end{pmatrix}\\
\ldots &\ldots &\ldots &\ldots &\ldots &\ldots \\
0 &0 & \ldots & 1 &\begin{pmatrix} t_{m-1}\\ t_m \end{pmatrix}& \begin{pmatrix} t_{m-1}\\ 2 \end{pmatrix}\\
0 &0 &\ldots & 0 & 1 & \begin{pmatrix} t_m\\ 2 \end{pmatrix}\end{matrix}\right|.
\end{equation}

Furthermore,
\newpage
$$
\left\{\begin{matrix} n\\k +2\end{matrix}\right\}=\left|\begin{matrix}\begin{pmatrix} n\\ t_1 \end{pmatrix} &
\begin{pmatrix} n\\ t_2 \end{pmatrix} & \ldots & \begin{pmatrix} n\\ t_{m-1} \end{pmatrix} &
\begin{pmatrix} n\\ t_m \end{pmatrix} & \begin{pmatrix} n\\ 2 \end{pmatrix} & 1\\
1 & \begin{pmatrix} t_1\\ t_2 \end{pmatrix}& \ldots & \begin{pmatrix} t_1\\ t_{m-1} \end{pmatrix}&
\begin{pmatrix} t_1\\ t_m\end{pmatrix}& \begin{pmatrix} t_1\\2 \end{pmatrix} & 1\\
0 & 1 & \ldots & \begin{pmatrix} t_2\\ t_{m-1} \end{pmatrix} & \begin{pmatrix} t_2\\t_m \end{pmatrix}&
\begin{pmatrix} t_2\\ 2 \end{pmatrix} & 1  \\
\ldots &\ldots & \ldots & \ldots & \ldots & \ldots & \ldots\\
0 & 0 & \ldots & 1 &\begin{pmatrix} t_{m-1}\\ t_m \end{pmatrix} &  \begin{pmatrix} t_{m-1}\\ 2 \end{pmatrix}&
1 \\
0 & 0 & \ldots &0& 1 &  \begin{pmatrix} t_m\\ 2 \end{pmatrix}&
1 \\
0 & 0 & \ldots & 0 & 0 & 1 & 1
\end{matrix}\right|=
$$

\begin{equation}\label{14}
-\left\{\begin{matrix} n\\k\end{matrix}\right\}+\left|\begin{matrix}\begin{pmatrix} n\\ t_1 \end{pmatrix}&
\begin{pmatrix} n\\ t_2 \end{pmatrix}&\ldots &\begin{pmatrix} n\\ t_{m-1} \end{pmatrix}&\begin{pmatrix} n\\ t_m \end{pmatrix}&\begin{pmatrix} n\\ 2 \end{pmatrix}\\
1& \begin{pmatrix} t_1\\ t_2 \end{pmatrix}& \ldots & \begin{pmatrix} t_1\\ t_{m-1} \end{pmatrix}&
\begin{pmatrix} t_1\\ t_m \end{pmatrix}& \begin{pmatrix} t_1\\ 2\end{pmatrix}\\
0 & 1 & \ldots & \begin{pmatrix} t_2\\ t_{m-1} \end{pmatrix}&\begin{pmatrix} t_2\\ t_m \end{pmatrix}&
\begin{pmatrix} t_2\\ 2\end{pmatrix}\\
\ldots & \ldots & \ldots & \ldots & \ldots & \ldots\\
0 & 0 & \ldots & 1 & \begin{pmatrix} t_{m-1}\\t_m \end{pmatrix}& \begin{pmatrix} t_{m-1}\\ 2 \end{pmatrix}\\
0 & 0 & \ldots & 0 & 1 & \begin{pmatrix} t_m\\ 2 \end{pmatrix}
\end{matrix}\right|.
\end{equation}

Thus, by (\ref{13}) and (\ref{14}), we have

$$
\left\{\begin{matrix} n\\k+3 \end{matrix}\right\}=\left\{\begin{matrix} n\\k \end{matrix}\right\}-\left|\begin{matrix}\begin{pmatrix} n\\ t_1 \end{pmatrix}&
\begin{pmatrix} n\\ t_2 \end{pmatrix}&\ldots &\begin{pmatrix} n\\ t_{m-1} \end{pmatrix}&\begin{pmatrix} n\\ t_m \end{pmatrix}&\begin{pmatrix} n\\ 1\end{pmatrix}\\
1& \begin{pmatrix} t_1\\ t_2 \end{pmatrix}& \ldots & \begin{pmatrix} t_1\\ t_{m-1} \end{pmatrix}&
\begin{pmatrix} t_1\\ t_m \end{pmatrix}& \begin{pmatrix} t_1\\ 1\end{pmatrix}\\
0 & 1 & \ldots & \begin{pmatrix} t_2\\ t_{m-1} \end{pmatrix}&\begin{pmatrix} t_2\\ t_m \end{pmatrix}&
\begin{pmatrix} t_2\\1\end{pmatrix}\\
\ldots & \ldots & \ldots & \ldots & \ldots & \ldots\\
0 & 0  &\ldots & 1 & \begin{pmatrix} t_{m-1}\\t_m \end{pmatrix}& \begin{pmatrix} t_{m-1}\\ 1 \end{pmatrix}\\
0 & 0 & \ldots & 0 & 1 & \begin{pmatrix} t_m\\ 1 \end{pmatrix}
\end{matrix}\right|+
$$
\newpage
\begin{equation}\label{15}
\left|\begin{matrix}\begin{pmatrix} n\\ t_1 \end{pmatrix}&
\begin{pmatrix} n\\ t_2 \end{pmatrix}&\ldots &\begin{pmatrix} n\\ t_{m-1} \end{pmatrix}&\begin{pmatrix} n\\ t_m \end{pmatrix}&\begin{pmatrix} n\\ 2\end{pmatrix}\\
1& \begin{pmatrix} t_1\\ t_2 \end{pmatrix}& \ldots & \begin{pmatrix} t_1\\ t_{m-1} \end{pmatrix}&
\begin{pmatrix} t_1\\ t_m \end{pmatrix}& \begin{pmatrix} t_1\\2\end{pmatrix}\\
0 & 1 & \ldots & \begin{pmatrix} t_2\\ t_{m-1} \end{pmatrix}&\begin{pmatrix} t_2\\ t_m \end{pmatrix}&
\begin{pmatrix} t_2\\2\end{pmatrix}\\
\ldots & \ldots & \ldots & \ldots & \ldots & \ldots\\
0 & 0  &\ldots & 1 & \begin{pmatrix} t_{m-1}\\t_m \end{pmatrix}& \begin{pmatrix} t_{m-1}\\ 2 \end{pmatrix}\\
0 & 0 & \ldots & 0 & 1 & \begin{pmatrix} t_m\\ 2 \end{pmatrix}
\end{matrix}\right|.
\end{equation}

Finally,
$$
\left\{\begin{matrix} n\\k+1 \end{matrix}\right\}=\left|\begin{matrix}\begin{pmatrix} n\\ t_1 \end{pmatrix}&
\begin{pmatrix} n\\ t_2 \end{pmatrix}&\ldots &\begin{pmatrix} n\\ t_{m-1} \end{pmatrix}&\begin{pmatrix} n\\ t_m \end{pmatrix}&\begin{pmatrix} n\\ 1\end{pmatrix}& 1\\
1& \begin{pmatrix} t_1\\ t_2 \end{pmatrix}& \ldots & \begin{pmatrix} t_1\\ t_{m-1} \end{pmatrix}&
\begin{pmatrix} t_1\\ t_m \end{pmatrix}& \begin{pmatrix} t_1\\ 1\end{pmatrix}& 1\\
0 & 1 & \ldots & \begin{pmatrix} t_2\\ t_{m-1} \end{pmatrix}&\begin{pmatrix} t_2\\ t_m \end{pmatrix}&
\begin{pmatrix} t_2\\1\end{pmatrix}& 1\\
\ldots & \ldots & \ldots & \ldots & \ldots & \ldots &\ldots\\
0 & 0  &\ldots & 0 &1 & \begin{pmatrix} t_m\\1\end{pmatrix}& 1 \\
0 & 0 & \ldots & 0 & 0 & 1 & 1
\end{matrix}\right|=
$$

\begin{equation}\label{16}
-\left\{\begin{matrix} n\\k \end{matrix}\right\}+\left|\begin{matrix}\begin{pmatrix} n\\ t_1 \end{pmatrix}&
\begin{pmatrix} n\\ t_2 \end{pmatrix}&\ldots &\begin{pmatrix} n\\ t_{m-1} \end{pmatrix}&\begin{pmatrix} n\\ t_m \end{pmatrix}&\begin{pmatrix} n\\ 1\end{pmatrix}\\
1& \begin{pmatrix} t_1\\ t_2 \end{pmatrix}& \ldots & \begin{pmatrix} t_1\\ t_{m-1} \end{pmatrix}&
\begin{pmatrix} t_1\\ t_m \end{pmatrix}& \begin{pmatrix} t_1\\ 1\end{pmatrix}\\
0 & 1 & \ldots & \begin{pmatrix} t_2\\ t_{m-1} \end{pmatrix}&\begin{pmatrix} t_2\\ t_m \end{pmatrix}&
\begin{pmatrix} t_2\\1\end{pmatrix}\\
\ldots & \ldots & \ldots & \ldots & \ldots & \ldots \\
0 & 0  &\ldots & 0 &1 & \begin{pmatrix} t_m\\1\end{pmatrix} \\
0 & 0 & \ldots & 0 & 0 & 1
\end{matrix}\right|.
\end{equation}

Subtracting (\ref{16}) from (\ref{14}) and, after that, subtracting the result from (\ref{15}), we find

$$
\left\{\begin{matrix} n\\k+3 \end{matrix}\right\}-\left(\left\{\begin{matrix} n\\k+2 \end{matrix}\right\}-
\left\{\begin{matrix} n\\k +1\end{matrix}\right\}\right)=\left\{\begin{matrix} n\\k \end{matrix}\right\},
$$

which completes proof. $\blacksquare $

Note that, in particular, we have\newpage

$$
\sum^{2^r-1}_{k=0}(-1)^k\left\{\begin{matrix} n\\k \end{matrix}\right\}=0, \;\; r\geq 2.
$$

Taking into account that (\cite{10})

$$
\sum^{2^r-1}_{k=0}\left\{\begin{matrix} n\\k \end{matrix}\right\}=n(n-1) \ldots (n-r+1), \;\;1\leq r \leq n-1,
$$

we find

\begin{equation}\label{17}
\sum^{2^r-1}_{i=0}\left\{\begin{matrix} n\\2i \end{matrix}\right\}=\sum^{2^r-1}_{i=0}\left\{\begin{matrix} n\\2i+1 \end{matrix}\right\}=\frac 1 2 n(n-1)\ldots(n-r+1),\;\; 2\leq r \leq n-1,
\end{equation}

as well and, in particular,

\begin{equation}\label{18}
\sum^{2^{n-1}-1}_{i=0}\left\{\begin{matrix} n\\2i \end{matrix}\right\}=
\sum^{2^{n-1}-1}_{i=0}\left\{\begin{matrix} n\\2i +1\end{matrix}\right\}=\frac {n!}{2}, \;\; n\geq 2.
\end{equation}

\section{Main conjecture}

     Denote by $\left[\begin{matrix} n\\k \end{matrix}\right]$ the number of all full cycles which have up-down index $k$. Many observations show that $\left[\begin{matrix} n\\k \end{matrix}\right]\approx\frac 1 n \left\{\begin{matrix} n\\k \end{matrix}\right\}$ with highly good approximation.Moreover, we think that the following conjecture is true.

\bfseries Conjecture 2. \mdseries     \slshape Let $t_1=t_1(k)$ be defined by (\ref{6}). If all the divisors of $n\geq 3$,
which are different from 1 are larger than $t_1,$ then exactly

\begin{equation}\label{19}
\left[\begin{matrix} n\\k \end{matrix}\right]=\frac 1 n \left(\left\{\begin{matrix} n\\k \end{matrix}\right\}-\left\{\begin{matrix} 0\\k \end{matrix}\right\}\right),
\end{equation}

where according to (\ref{10}), $\left\{\begin{matrix} 0\\k \end{matrix}\right\}=\tau_k;$\newline
\indent otherwise, (19) is a good approximation of $\left[\begin{matrix} n\\k \end{matrix}\right].$\upshape \newline
\indent Note that, in the conditions of Conjecture 2 the fraction in  (\ref{19}) is an integer. Indeed,
from  (\ref{11}) we find
\newpage
\begin{equation}\label{20}
\left\{\begin{matrix} n\\k\end{matrix}\right\}- \left\{\begin{matrix} 0\\k\end{matrix}\right\}=
\left|\begin{matrix}
\begin{pmatrix} n\\ t_1 \end{pmatrix}&\begin{pmatrix} n\\ t_2 \end{pmatrix}&\begin{pmatrix} n\\ t_3 \end{pmatrix}&\ldots &\begin{pmatrix} n\\ t_{m-1} \end{pmatrix}&\begin{pmatrix} n\\ t_m \end{pmatrix}& 0\\
1& \begin{pmatrix} t_1\\ t_2 \end{pmatrix}& \begin{pmatrix} t_1\\ t_3 \end{pmatrix}&\ldots & \begin{pmatrix} t_1\\ t_{m-1} \end{pmatrix}&\begin{pmatrix} t_1\\ t_m \end{pmatrix}&  1\\
0 & 1 & \begin{pmatrix} t_2\\ t_3 \end{pmatrix}&\ldots & \begin{pmatrix} t_2\\ t_{m-1} \end{pmatrix}&\begin{pmatrix} t_2\\ t_m \end{pmatrix}&1\\
0 & 0 & 1 & \ldots &\begin{pmatrix} t_3\\ t_{m-1} \end{pmatrix}&\begin{pmatrix} t_3\\ t_m \end{pmatrix}& 1\\
\ldots & \ldots & \ldots & \ldots & \ldots & \ldots &\ldots\\
0 & 0 &0 &\ldots & 1 &\begin{pmatrix} t_{m-1}\\ t_m \end{pmatrix} & 1
\\0 & 0 &0& \ldots & 0 & 1 & 1
\end{matrix}\right|.
\end{equation}

Matrix (\ref{20}) differs from matrix (\ref{11}) only in the last element of the first row. To this element corresponds the only diagonal in matrix (\ref{11}) composed of 1's. The corresponding term in determinant  (\ref{11}) is

$$
(-1)^m=\tau_k=\left\{\begin{matrix} 0\\k\end{matrix}\right\}.
$$

In the conditions of Conjecture 2 all elements of the first row of matrix (\ref{20}) are divided by $n$. Therefore, in this case, the fraction in (\ref{19}) is an integer. These arguments allow to prove the following statement.
\begin{theorem}\label{7}
If Conjecture $2$ is true, then, in contrast of sequence $(\left\{\begin{matrix} n\\k\end{matrix}\right\})_{n\geq 0},$ for some values of $k$ the sequence $ (\left[\begin{matrix} n\\k\end{matrix}\right])_{n\geq 0}$
is not polynomial.
 \end{theorem}
 \bfseries Proof. \mdseries Indeed, let sequence $(\left[\begin{matrix} n\\k\end{matrix}\right])$ be polynomial for every $k:$ $\left[\begin{matrix} n\\k\end{matrix}\right]=P^{(k)}(n).$  Denoting $Q^{(k)}(n)$  polynomial (20) and considering prime values of $n,$ we conclude that $P^{(k)}(n)\equiv Q^{(k})(n)/n.$ Besides, since $\left[\begin{matrix} n\\k\end{matrix}\right]$ is integer for every $n,$ then $ Q^{(k)}(n)/n$ is integer-valued polynomial in $n$ for every $k.$ Let $k$ be of the form $k=2^{p-1},$ where $p$ is prime. Then $t_1=p,$ and in(11) $m=1,$ i.e., we have a matrix $2\times2:$
 $$\left\{\begin{matrix} n\\
2^{p-1}\end{matrix}\right\}=\begin{vmatrix}  \begin{pmatrix}n\\p\end{pmatrix}&1\\
1 & 1\end{vmatrix}= \begin{pmatrix}n\\p\end{pmatrix}-1$$
and, by (20),
 $$Q(n)/n=\frac {1} {n}\begin{vmatrix}  \begin{pmatrix}n\\p\end{pmatrix}&0\\
1 & 1\end{vmatrix}= \frac{(n-1)(n-2)\ldots (n-p+1)}{p!}\enskip.$$\newpage
Let $n$ be multiple of p. We see that, for such $n,\enskip Q(n)/n$ is not integer. We have a contradiction which completes proof. $\blacksquare$

     Furthermore, in connection with Theorem 6, note that, if Conjecture 2 is true, then in its conditions for
     $ k \equiv 0 \pmod 4$ we also have

 \begin{equation}\label{21}
\left[\begin{matrix} n\\k \end{matrix}\right] -\left[\begin{matrix} n\\k+1 \end{matrix}\right]
+\left[\begin{matrix} n\\k+2 \end{matrix}\right]-\left[\begin{matrix} n\\k+3 \end{matrix}\right]=0,\;\; n\geq 3.
\end{equation}

In particular, (\ref{21}) is true for $n$ being an odd prime.

\section{An analog of sequence A360651 for full cycles}

Put for $k\geq 1$

\begin{equation}\label{22}
g(k)=\lfloor\log_2 k\rfloor+1,\;\;\;h(k)=k-2^{g(k)-1}.
\end{equation}

Sequence A360051 \cite{11} is the sequence

$$
\left\{\begin{matrix} 1\\0\end{matrix}\right\}, \left\{\begin{matrix} 2\\0\end{matrix}\right\},
\left\{\begin{matrix} 2\\1\end{matrix}\right\},\left\{\begin{matrix} 3\\0\end{matrix}\right\},
\left\{\begin{matrix} 3\\1\end{matrix}\right\},\left\{\begin{matrix} 3\\2\end{matrix}\right\},
\left\{\begin{matrix} 3\\3\end{matrix}\right\},
$$

\begin{equation}\label{23}
\left\{\begin{matrix}4\\0\end{matrix}\right\},\left\{\begin{matrix} 4\\1\end{matrix}\right\},
\left\{\begin{matrix} 4\\2\end{matrix}\right\}, \left\{\begin{matrix} 4\\3\end{matrix}\right\},
\left\{\begin{matrix} 4\\4\end{matrix}\right\},\left\{\begin{matrix} 4\\5\end{matrix}\right\},
\left\{\begin{matrix} 4\\6\end{matrix}\right\},\left\{\begin{matrix} 4\\7\end{matrix}\right\},\ldots ,
\end{equation}

i.e. the sequence \cite{10}

\begin{equation}\label{24}
\left(\left\{\begin{matrix} g(k)\\h(k)\end{matrix}\right\}\right)_{k\geq1}.
\end{equation}

     The most simple algorithm for evaluation of $\left\{\begin{matrix} n\\k\end{matrix}\right\}$
is the following recursion which is directly obtained from Theorem 16 \cite{10}.

\begin{theorem}\label{t8}   We have

\begin{equation}\label{25}
\left\{\begin{matrix} n\\k\end{matrix}\right\}= \left\{\begin{matrix} g(k)\\h(k)\end{matrix}\right\}\begin{pmatrix}
n\\g(k)\end{pmatrix}-\left\{\begin{matrix} n\\h(k)\end{matrix}\right\}, \;\; k\geq 1
\end{equation}
\end{theorem}

with the  initial condition $\left\{\begin{matrix} n\\0\end{matrix}\right\}=1$.

     An analog of sequence A360051 for full cycles is the sequence
$$
\left[\begin{matrix} 1\\0\end{matrix}\right],  \left[\begin{matrix} 2\\0\end{matrix}\right],    \left[\begin{matrix} 2\\1\end{matrix}\right],\left[\begin{matrix} 3\\0\end{matrix}\right],
\left[\begin{matrix} 3\\1\end{matrix}\right],  \left[\begin{matrix} 3\\2\end{matrix}\right],   \left[\begin{matrix} 3\\3\end{matrix}\right],
$$

\begin{equation}\label{26}
\left[\begin{matrix} 4\\0\end{matrix}\right],\left[\begin{matrix} 4\\1\end{matrix}\right],
\left[\begin{matrix} 4\\2\end{matrix}\right],\left[\begin{matrix} 4\\3\end{matrix}\right],
\left[\begin{matrix} 4\\4\end{matrix}\right],\left[\begin{matrix} 4\\5\end{matrix}\right],
\left[\begin{matrix} 4\\6\end{matrix}\right],\left[\begin{matrix} 4\\7\end{matrix}\right],
\hdots,
\end{equation}
\newpage
i.e., the sequence

\begin{equation}\label{27}
\left(\left[\begin{matrix} g(k)\\h(k)\end{matrix}\right]\right)^{\infty}_{k=1}.
\end{equation}

 By convention, $\left[\begin{matrix} 1\\0\end{matrix}\right]=1$. Note that, for $n=2,$ there is only cycle $(2,1)$ in case $k=0$,
i.e. $\left[\begin{matrix} 2\\0\end{matrix}\right]=1,\;\;\left[\begin{matrix} 2\\1\end{matrix}\right]=0$.
Note that, for $n\geq 3,\;\;\left[\begin{matrix} n\\0\end{matrix}\right]=0$, since the only permutation corresponding to this case has more than one cycle. Furthermore, similar to $\left\{\begin{matrix} n\\k\end{matrix}\right\}$ \cite{10}, one can prove that, for $n\geq 3$

\begin{equation}\label{28}
\left[\begin{matrix} n\\k\end{matrix}\right]=\left[\begin{matrix} n\\2^{n-1}-1-k\end{matrix}\right],\;\; 0\leq k\leq 2^{n-2}-1.
\end{equation}

In particular, each block of sequence (\ref{26}) begins and ends with 0:

\begin{equation}\label{29}
\left[\begin{matrix} n\\2^{n-1}-1\end{matrix}\right]=\left[\begin{matrix} n\\0\end{matrix}\right]=0,\;\; n\geq 3.
\end{equation}

Note that the conditions of Conjecture 2 are satisfied for a \slshape whole \upshape block
$\left(\left[\begin{matrix} n\\k\end{matrix}\right]\right)^{2^{n-1}-1}_{k=0}$ if and only if $n$ is an odd prime. For example, for $n=3$ we have only two full cycles $(2,3,1)$ and $(3,1,2$ with $k=2$ and $k=1$ correspondingly. Thus, this block in (\ref{26}) has the form: $0,1,1,0$.  This conforms to (\ref{20}). Indeed,

$$
\left[\begin{matrix} 3\\1\end{matrix}\right]=\frac 1 3 \left|\begin{matrix}
\begin{pmatrix} 3\\1 \end{pmatrix}&0\\ 1& 1
\end{matrix}\right|=1,\;\;\left[\begin{matrix} 3\\2\end{matrix}\right]=\frac 1 3 \left|\begin{matrix}
\begin{pmatrix} 3\\2 \end{pmatrix}&0\\ 1& 1
\end{matrix}\right|=1,
$$

$$
\left[\begin{matrix} 3\\3\end{matrix}\right]=\frac 1 3 \left|\begin{matrix}
\begin{pmatrix} 3\\2 \end{pmatrix}&\begin{pmatrix} 3\\1 \end{pmatrix}&0\\ 1&\begin{pmatrix} 2\\1 \end{pmatrix}&1\\
0 & 1 & 1
\end{matrix}\right|=0.
$$

For $n=5$ we have the only cycle $(5,4,2,1,3)$ with $k=1$, two cycles $(4,3,1,5,2)$ and $(5,4,1,3,2)$ with $k=2$,
only cycle $(5,3,1,2,4)$ with $k=3$, two cycles $(4,1,5,3,2)$ and $(4,3,5,2,1)$ with $k=4$, three cycles $(3,1,5,2,4), (5,1,4,2,3)$ and $(5,3,4,2,1)$ with $k=5$, two cycles $(3,1,4,5,2)$ and $(4,1,2,5,3)$ with $k=6$, only cycle $(5,1,2,3,4)$ with $k=7$ and the numbers of cycles with $k=8,9,\ldots , 14$ conform to (\ref{28}). They are:\newpage

$$
(3,5,4,2,1);\;(2,5,4,1,2),\;(3,5,2,1,4);\;(2,4,1,5,3),\;(3,4,2,5,1),
$$
$$(4,5,2,3,1); \; (2,5,1,3,4),\;(4,5,1,2,3);\;
(2,4,5,3,1); \;(3,4,5,1,2),
$$
$$(2,3,5,1,4);\;(2,3,4,5,1)$$.

Thus, this block in (\ref{26}) has the form:
$0,1,2,1,2,3,2,1,1,2,3,2,1,2,1,0$. This conforms to (\ref{19}), (\ref{20}) and (\ref{21}). Indeed,

$$
\left[\begin{matrix}5\\1\end{matrix}\right]= \frac 1 5\left|\begin{matrix}
\begin{pmatrix} 5\\1 \end{pmatrix}&0\\ 1&1\\
\end{matrix}\right|=1,\;\;\left[\begin{matrix}5\\2\end{matrix}\right]= \frac 1 5\left|\begin{matrix}
\begin{pmatrix} 5\\2 \end{pmatrix}&0\\ 1&1\\
\end{matrix}\right|=2,
$$

$$
\left[\begin{matrix}5\\3\end{matrix}\right]= \frac 1 5\left|\begin{matrix}
\begin{pmatrix} 5\\2 \end{pmatrix}&\begin{pmatrix} 5\\1 \end{pmatrix}&0\\ 1&\begin{pmatrix} 2\\1 \end{pmatrix}&1\\
0&1 & 1
\end{matrix}\right|=1,\;\; \left[\begin{matrix}5\\4\end{matrix}\right]= \frac 1 5\left|\begin{matrix}
\begin{pmatrix} 5\\3 \end{pmatrix}&0\\ 1&1\\
\end{matrix}\right|=2,
$$

$$
\left[\begin{matrix}5\\5\end{matrix}\right]= \frac 1 5\left|\begin{matrix}
\begin{pmatrix} 5\\3 \end{pmatrix}&\begin{pmatrix} 5\\1 \end{pmatrix}&0\\ 1&\begin{pmatrix}3\\1 \end{pmatrix}&1\\
0&1 & 1
\end{matrix}\right|=3,\;\;\left[\begin{matrix}5\\6\end{matrix}\right]= \frac 1 5\left|\begin{matrix}
\begin{pmatrix} 5\\3 \end{pmatrix}&\begin{pmatrix} 5\\2 \end{pmatrix}&0\\ 1&\begin{pmatrix}3\\2 \end{pmatrix}&1\\
0&1 & 1
\end{matrix}\right|=2,
$$

$$
\left[\begin{matrix}5\\7\end{matrix}\right]= \frac 1 5\left|\begin{matrix}
\begin{pmatrix} 5\\3 \end{pmatrix}&\begin{pmatrix} 5\\2 \end{pmatrix}&\begin{pmatrix} 5\\1 \end{pmatrix}&0\\ 1&\begin{pmatrix} 3\\2 \end{pmatrix}&\begin{pmatrix}3\\1 \end{pmatrix}&1\\
0&1&\begin{pmatrix} 2\\1 \end{pmatrix} & 1\\
0 & 0 & 1 & 1
\end{matrix}\right|=1,\;\; etc.
$$

At the same time, for $n=4$ which is not a prime, (\ref{20}) is satisfied, generally speaking, only approximately.
Indeed, here we have only permutations with indices
$$k=1,2, \ldots ,6:\; (4,3,1,2), (3,1,4,2), (4,1,2,3)
(3,4,2,1), (2,4,1,3), (2,3,4,1)
$$
correspondingly. I.e. this block in (\ref{26}) hast the form $0,1,1,1,1,1,1,0$, while according to (\ref{20}) we have
$0,1,\frac 3 2, \frac 1 2, 1,1, \frac 1 2, \frac 1 2$.

This, the first numbers of sequence (\ref{26}) are:

\begin{equation}\label{30}
1,1,0,0,1,1,0,0,1,1,1,1,1,1,0,0,1,2,1,2,3,2,1,1,2,3,2,1,2,1,0,0,\ldots
\end{equation}
\newpage

\section{Some other open problems}

1. It is very interesting to estimate the remainder term of approximation (\ref{19}) in the general case.

2. Let in the block $"n"(n\geq 3)$ in (\ref{26}) for every $k\in[0,2^{n-1}-1]$ which is divided by 4, (\ref{21})
   be satisfied. We conjecture that in this case $n$ is a prime.

3. It is known that the sequence $(a_n)$ of the numbers of the alternating permutations of elements $1,2,\ldots, n$
   for which $\pi_1< \pi_2 > \pi_3 < \ldots$ for $n \geq 1$ is (cf. A000111[11]),

\begin{equation}\label{31}
1,1,2,5,16,61,272,1385,7936,50521,353792,2702765,\ldots
\end{equation}

     The corresponding sequence $(f_n)$ of the numbers of alternating full cycles is (\cite{3}):

\begin{equation}\label{32}
1,0,1,1,3,10,39,173,882,5052,32163,225230,\ldots
\end{equation}

It is naturally to conjecture (\cite{3}, 1996) that

\begin{equation}\label{33}
f_n\approx \frac{a_n}{n}.
\end{equation}

Indeed, the sequence $\left(\frac{a_n}{n}\right)$ gives a highly good approximation of (\ref{32}):

\begin{equation}\label{34}
1,\;0.5,\;0.7,\;1.3,\;3.3,\;10.2,\;38.9,\;173.1,\;881.8,\;5052.1,\;32162.9,\;225230.4,\; \ldots.
\end{equation}

It is interesting to prove (\ref{33}) with an estimate of the remainder term.

4. A difficult combinatorial problem - to enumerate the alternating permutations and \slshape antialternating permutations \upshape for which $\pi_1> \pi_2 < \pi_3 > \ldots$ without fixed points or, the same, without cycles of length 1. The first numbers of these sequences $(b_n)$ and $(b_n^*)$ for $n \geq 1$ are

\begin{equation}\label{35}
0,\;0,\;1,\;2,\;6,\;22,\;102,\;506,\;2952,\;18502,\;131112,\;991226,\;\ldots
\end{equation}

and

\begin{equation}\label{36}
0,\;1,\;1,\;2,\;6,\;24,\;102,\;528,\;2952,\;19008,\;131112,\;1009728,\;\ldots
\end{equation}

It is not difficult to prove that

\begin{equation}\label{37}
b_{2n-1}^*=b_{2n-1}, \;\; n\geq 1,
\end{equation}\newpage

\begin{equation}\label{38}
b_{2n}^*=b_{2n}+b_{2n-2}, \;\; n\geq 2.
\end{equation}
i.e., for permutations without fixed points, we have a violation of the symmetry.

\bfseries Conjecture 3.\mdseries

\begin{equation}\label{39}
\lim_{n\rightarrow\infty}\frac{a_n}{b_n}=\lim_{n\rightarrow\infty}\frac{a_n}{b_n^*}=e.
\end{equation}

Moreover, we conjecture that, since $n=3$

\begin{equation}\label{40}
\frac{a_{2n}}{b^*_{2n}}< e < \frac{a_{2n}}{b_{2n}},
\end{equation}

such that $\frac{a_{2n}}{b_{2n}^*}$ increases and $\frac{a_{2n}}{b_{2n}}$ decreases. So, for $n=3,$

$$
\frac{61}{24}=2.541\hdots < e < \frac{61}{22}= 2.772 \hdots;
$$

for $n=4,$

$$
\frac{1385}{528}=2.623 \hdots < e < \frac{1385}{506}=2.737 \hdots;
$$

for $n=5,$

$$
\frac{50521}{19008}=2.657\hdots < e < \frac{50521}{18502}= 2.730 \hdots;
$$

for $n=6,$

$$
\frac{2702765}{1009728}=2.676\hdots < e < \frac{2702765}{991226}= 2.726\hdots\;\;\; etc.
$$

If to consider the concatenation sequence similar to (\ref{26}) for permutations having no fixed points with up-down index $k\geq 0,$ then we obtain a sequence which is asymmetric in its blocks (except $n=3$) with the following first numbers:

$$
0,\;1,\;0,\;0,\;1,\;1,\;0,\;1,\;1,\;2,\;1,\;1,\;2,\;1,\;0,\;0,\;2,
$$
\begin{equation}\label{41}
4,\;2,\;5,\;5,\;4,\;1,\;2,\;4,\;6,\;3,\;2,\;3,\;1,\;0,\;\hdots
\end{equation}

5. An algorithm for calculation of the cyclic indicators for the alternating and antialternating  permutations with restricted positions, given by any (0,1) matrix A, was obtained in \cite{3} with its realization in Turbo-Pascal 6.0.
For example, if $A=J_6$ (i.e. without restrictions on positions) the "alternating" indicator has the form

$$
10t_6+12t_1t_5+10t_1t_2t_3+8t_1^2t_4+7t_2t_4+4t_1^2t_2^2+4t_3^2+4t_1^3t_3+t_2^3+t_1^4t_2,
$$
\newpage
while the "antialternating" indicator has the form

$$
10t_6+12t_1t_5+12t_1t_2t_3+7t_1^2t_4+8t_2t_4+4t_1^2t_2^2+4t_3^2+2t_1^3t_3+2t_2^3,
$$

such that the difference between these indicators is

\begin{equation}\label{42}
R_6=(t_1^2-t_2)(t_4+2t_1t_3+t_2(t_1^2+t_2)).
\end{equation}

Analogously, we have

\begin{equation}\label{43}
R_2=t_1^2-t_2,\;\;R_4=0,
\end{equation}

\begin{equation}\label{44}
R_8=(t_1^2-t_2)(10t_6+12t_1t_8+8t_1t_2t_3+4t_3^2+2t_1^2t_2^2+(t_1^2+t_2)(7t_4+2t_1t_3+t_2^2))
\end{equation}
$$
R_{10}=(t_1^2-t_2)(173t_8+198t_1t_7+120t_1t_3t_4+96t_3t_5+
$$
$$
+96t_1t_2t_5+43t_4^2+39t_1^2t_2t_4+26t_1t_2^2t_3+9t_1^2t_2^3+
$$
$$
+(t_1^2+t_2)(110t_6+40t_3^2+36t_1t_5+34t_2t_4+30t_1t_2t_3+t_2^3)+
$$
\begin{equation}\label{45}
+(t_1^4+t_1^2t_2+t_2^2)(3t_2^2+6t_4)), etc.
\end{equation}

     In addition, it is easy to see that, for odd $n$, we have $R_n=0$. In case $t_1=0,\;t_2=\hdots=t_{2n}=1,$  we
obtain that $R_{2n}=b_{2n-2}$ (sf.(\ref{35}) for $n\geq 2$). It follows from (\ref{38}).

\bfseries Conjecture 4.\mdseries  \slshape \enskip Polynomial $R_{2n}$ is multiple of $t_1^2-t_2$ and, moreover, all coefficients of the polynomial $\frac{R_{2n}}{t_1^2-t_2}$ are positive.\upshape

Note that, the sequence of the maximal coefficients of polynomials $\frac{R_{2n}}{t_1^2-t_2},\;\;n=1,2,\hdots$,
is

\begin{equation}\label{46}
1,0,2,12,198,\hdots
\end{equation}

Whether the maximal coefficient of the polynomial $\frac{R_{2n}}{t_1^2-t_2}$ is always the coefficient of $t_1t_{2n-3}?$

Note that if Conjecture 4 is true, then \slshape the numbers of the alternating and the antialternating full cycles are equal for $n\geq 3$.\upshape  Indeed, if $t_1=t_2=0$ then, by Conjecture 4, always $R_n=0$.

Finally, if Conjecture 4 is true, then for $t_1^2=t_2=t^2$ we have

\begin{equation}\label{47}
R_{2n}(t,t^2,t_3,t_4, \hdots)=0.
\end{equation}
\newpage
The latter means that \slshape the numbers of all alternating and all antialternating permutations of elements $1,2,\hdots,2n$
having the same given summary length of cycles of length 1 and 2 and the same given numbers of cycles of length
$i,\;\; i=3,4,\hdots, 2n$, are equal.\upshape

     E.g., in case of $n=10,$ the sum of the coefficients of $t_1^6t_4,\;t_1^4t_2t_4$, $t_1^2t_2^2t_4$
and $t_2^3t_4$ in the
"alternating" indicator is
$$
6+241+770+248=1265
$$
and in the "antialternating" indicator it is

$$
0+168+809+288=1265.
$$

\end {document}